\documentclass[12pt]{amsart}
\usepackage{amsfonts, amssymb,amsmath, amsthm, amsxtra, latexsym, amscd}
\usepackage[latin1]{inputenc} 
\usepackage{color}

\def\draft{\centerline{(Draft {\the \day}/{\the\month} \the \year.)}}

\theoremstyle{definition}
\newtheorem{theo+}    {Theorem}      [section]
\newtheorem{prop+}  [theo+]  {Proposition}
\newtheorem{coro+}  [theo+]  {Corollary}
\newtheorem{lemm+}  [theo+]  {Lemma}
\newtheorem{deep+}  [theo+]  {Deep Result}
\newtheorem{fact+}  [theo+]  {Fact}
\theoremstyle{definition}
\newtheorem{exam+}  [theo+]  {Example}
\newtheorem{rema+}  [theo+]  {Remark}
\newtheorem{defi+}  [theo+]  {Definition}
\newtheorem{xca+}[theo+]{Exercise}

\numberwithin{equation}{section}

\def\draft{\centerline{(Draft {\the \day}/{\the\month} \the \year.)}}
\bigskip

\def\refn#1.#2{\expandafter\def\csname#1\endcsname{[#2]}}
\def\refnr#1.{\csname#1\endcsname}

\def\fa{\mathfrak a}

\def\fg{\mathfrak g}
\def\fk{\mathfrak k}

\def\fp{\mathfrak p}

\def\a{\alpha}

\def\Claminv2{|C(\Lambda)|^{-2}}

\def\varepsi{\varepsilon}
\def\lam{\lambda}

\def\de{d\varepsilon}

\def\Aa2D{A^{\a,2}(D)}
\def\bAa2D{\overline{A^{\a,2}(D)}}
\def\Ab2D{A^{\beta,2}(D)}
\def\bAb2D{\overline{A^{\beta,2}(D)}}

\def\Norm#1_#2{\Vert#1\Vert_{#2}}

\def\phipl12{\phi_{p_{l_1}, p_{l_2}}}
\def\phip01{\phi_{p_{0}, p_{0}}}

\def\a{\alpha}

\def\Claminv2{|C(\Lambda)|^{-2}}

\def\varepsi{\varepsilon}

\def\sig{\sigma}

\def\lam{\lambda}

\def\det{\operatorname{det}}
\def\diag{\operatorname{diag}}
\def\Ind{\operatorname{Ind}}
\def\diag{\operatorname{diag}}

\def\exp{\operatorname{exp}}

\def\Pr{\operatorname{Pr}}

\def\vol{\operatorname{vol}}

\def\de{d\varepsilon}

\def\Aa2D{A^{\a,2}(D)}
\def\bAa2D{\overline{A^{\a,2}(D)}}
\def\Ab2D{A^{\beta,2}(D)}
\def\bAb2D{\overline{A^{\beta,2}(D)}}

\def\phipl12{\phi_{p_{l_1}, p_{l_2}}}
\def\phip01{\phi_{p_{0}, p_{0}}}
\def\m{\underline{\bold m}}

\def\bc{\mathbb C}
\def\br{\mathbb R}

\def\alg/{algebra} 
\def\Alg/{Algebra} 
\def\alt/{alternative} 
\def\anal/{analytic}
\def\analfunc/{\anal/\ \func/}
\def\Ans/{\it Answer. \normal}
\def\ass/{associative}
\def\nass/{non-\ass/}
\def\autom/{automorphism}
\def\homom/{homomorphism}
\def\isom/{isomorphism}
\def\bdd/{bounded}
\def\Bdd/{Bounded}
\def\bddsymdom/{bounded \sym/ \dom/}
\def\Cartdom/{Cartan \dom/}
\def\bdry/{boundary}
\def\bsd/{\bdd/ \symdom/}
\def\bv/{boundary value}
\def\cf/{{\it cf}\.}
\def\Cf/{{\it Cf}\.}
\def\charr/{character}
\def\coeff/{coefficient}
\def\comm/{commutative}
\def\cpct/{compact}
\def\compl/{complex}
\def\comp/{complex}
\def\Comp/{Complex}
\def\conf/{conformal}
\def\conj/{conjugate}
\def\conn/{connect}
\def\cont/{continuous}
\def\conv/{converge} 
\def\convc/{convergence}
\def\convt/{convergent}
\def\convx/{convex}
\def\coord/{coordinate}
\def\lcoord/{local coordinate}
\def\Corr/{Corresponding}
\def\corr/{corresponding}
\def\corrd/{correspond}
\def\cov/{covariant}
\def\decomp/{decomposition}
\def\deco/{decompose}
\def\diff/{different} 
\def\Diff/{Different} 
\def\dimn/{dimension} 
\def\distr/{distribution} 
\def\div/{diverge} 
\def\dom/{domain}
\def\eg/{\hbox{\it e.g}\.}
\def\eigenf/{eigen\-\func/}
\def\eigensp/{eigen\-space}
\def\eigenv/{eigen\-value}
\def\eq/{equation}
\def\equa/{equation}
\def\de/{\diff/ial \equa/}
\def\do/{\diff/ial operator}
\def\ode/{ordinary \de/}
\def\pde/{partial \de/}
\def\pdo/{partial \diff/ial operator}
\def\psdo/{pseudo \diff/ial operator}
\def\fin/{finite}
\def\Ex/{\it Example.\ \normal}
\def\Exnr#1/{\it Example #1.\ \normal}
\def\foll/{follow}
\def\follg/{following}
\def\Follg/{Following}
\def\func/{function}
\def\Func/{Function}
\def\Fonc/{Fonc\-tion}
\def\fonc/{fonc\-tion}
\def\Funk/{Funk\-tion}
\def\funk/{Funk\-tion}
\def\gen/{general}
\def\har/{harmonic}
\def\Hint/{\it Hint. \normal}
\def\hist/{historic}
\def\histcl/{historical}
\def\hol/{holo\-morphic}
\def\homog/{ho\-mo\-ge\-ne\-ous}
\def\hyp/{hyper\-bolic}
\def\hyperg/{hyper\-geometric}
\def\ie/{\hbox{\it i.e.}}
\def\iff/{if and only if}
\def\ineq/{inequality}
\def\infra/{{\it inf\-ra}}
\def\ultra/{{\it ult\-ra}}
\def\Inpart/{In particular}
\def\inpart/{in particular}
\def\instof/{instead of}
\def\interps/{interpolation space}
\def\interp/{interpolation}
\def\Interp/{Interpolation}
\def\interpr/{Interpretation}
\def\Intr/{Introduction}
\def\intv/{interval}
\def\inv/{invariant}
\def\invc/{invariance}
\def\Iowords/{In other words}
\def\iowords/{in other words}
\def\ipr/{inner product}
\def\irred/{irreducible}
\def\lb/{line bundle}
\def\lin/{linear}
\def\lhs/{left hand side}
\def\rhs/{right hand side}
\def\loc/{local}
\def\math/{mathematic} 
\def\mathcn/{\math/ian}
\def\manif/{manifold}
\def\meas/{measure}
\def\measl/{measurable}
\def\mero/{mero\-morphic}
\def\mon/{monomial}
\def\monog/{monogenic}
\def\mult/{multiple}
\def\multy/{multiply}
\def\multn/{multiplication}
\def\nas/{necessary and sufficient}
\def\nbd/{neighborhood}
\def\neg/{negative}
\def\nondeg/{nondegenerate}
\def\Oohand/{On the other hand}
\def\oohand/{on the other hand}
\def\Oonhand/{On the one hand}
\def\oonhand/{on the one hand}
\def\oper/{operator}
\def\orth/{ortho\-gonal}
\def\orthon/{ortho\-normal}
\def\otoh/{on the other hand}
\def\quat/{quaternion}
\def\pp/{\hbox{a. e.}}
\def\psh/{plurisubharmonic}
\def\pol/{polynomial}
\def\pot/{potential}
\def\pos/{positive}
\def\princ/{principle}
\def\prob/{probability}
\def\proj/{projective}
\def\projn/{projection}
\def\Proof/{\it Proof:\normal}
\def\Rem/{\it Remark\normal}
\def\Remnr#1/{\it Remark\ \normal #1. }
\def\rep/{representation}
\def\meta/{metaplectic representation}
\def\repr/{reproducing}
\def\reprker/{reproducing kernel}
\def\resp/{respective} 
\def\resply/{respectively}
\def\restr/{restriction}
\def\sa/{self-adjoint}
\def\st/{such that}

\def\sol/{solution}
\def\ru/{space}
\def\sph/{spherical}
\def\ssp/{sub\ru/}
\def\sym/{symmetric}
\def\Sym/{Symmetric}
\def\symb/{symbol}
\def\symbc/{symbolic}
\def\symdom/{\sym/ domain}
\def\symp/{symplectic}
\def\Theor#1/{\fet Theorem #1.\ \normal}
\def\Lem#1/{\fet Lemma #1.\ \normal}
\def\Lemma/{\fet Lemma.\ \normal}
\def\topl/{topology}
\def\topll/{topological}
\def\transf/{transform}
\def\transl/{translation}
\def\transfn/{transformation}
\def\transv/{transvectant}
\def\trig/{trigonometric}
\def\tril/{trilinear}
\def\trilf/{trilinear form}
\def\uhp/{upper halfplane}
\def\uhs/{upper halfspace}
\def\vb/{vector bundle}
\def\vf/{vector field}
\def\vsp/{vector space}
\def\wrt/{with respect to}
\def\Wlog/{Without loss of generality}
\def\a{\alpha}

\def\lam{\lambda}

\def\sig{\sigma}

\def\Ab/{Abel}
\def\Ban/{Banach}
\def\Bansp/{\Ban/ space}
\def\Belt/{Bel\-tra\-mi}
\def\Berg/{Berg\-man}
\def\Bern/{Ber\-nou\-lli}
\def\Berz/{Berezin}
\def\Bess/{Bessel}
\def\Cart/{Car\-tan}
\def\Cay/{Cay\-ley}
\def\CG/{Clebsch-Gordan}
\def\Cl/{Clifford}
\def\CR/{Cauchy-Rie\-mann}
\def\Dir/{Dirichlet}
\def\Eucl/{Euclide}
\def\Eucln/{Euclidean}
\def\F/{Fourier}
\def\Hank/{Hankel}
\def\Hankf/{\Hank/ form}
\def\Herm/{Hermite}
\def\Hilb/{Hilbert}
\def\Hilbs/{Hilbert space}
\def\Hilbsp/{Hilbert space}
\def\HS/{Hilbert-Schmidt}
\def\Lag/{La\-grange}
\def\Lap/{La\-place}
\def\LapBelt/{\Lap/-\Belt/}
\def\Leb/{Lebesgue}
\def\Marc/{Mar\-cin\-kie\-wicz}
\def\Moeb/{Moebius}
\def\Moebt/{Moebius transformation}
\def\Moebtransfn/{Moebius transformation}
\def\Pla/{Plan\-che\-rel}
\def\Poin/{Poin\-car\'e}
\def\Riem/{Rie\-mann}
\def\Riemn/{\Riem/ian}
\def\psRiemn/{pseudo-\Riem/ian}
\def\Riems/{Rie\-mann surface}
\def\Schroe/{Schr\"odinger}
\def\Weier/{Weier\-strass}
\def\anal/{analytic}
\def\bsd/{bounded symmetric domain  }
\def\bdd/{bounded}
\def\calc/{calculation}\def\conj{conjugate}
\def\calci/{calculating}\def\eg{e.g.}
\def\conj/{conjugate}
\def\deco/{decomposition}
\def\eg/{e.g.}
\def\fct/{function}
\def\gp/{group}
\def\hw/{highest weight}
\def\hwv/{highest weight vector}
\def\hwvs/{highest weight vectors}
\def\lw/{lowest weight}
\def\lwv/{lowest weight vector}
\def\lwvs/{lowest weight vectors}
\def\hds/{holomorphic discrete series}
\def\iff/{if and only if}
\def\inv/{invariant}
\def\irrde/{irreducible decomposition}
\def\meas/{measure}
\def\transf/{transform}
\def\rep/{representation}
\def\resp/{respectively}
\def\inters/{intertwines}
\def\interg/{intertwining}
\def\meta/{metaplectic representation}
\def\qu/{quaternion}
\def\rep/{representation}
\def\symdom/{ symmetric domain}
\def\st/{such that}
\def\shd/{subhead}
\def\transf/{transform}
\def\wrt/{with respect to}

\def\Norm#1#2#3{\Vert#1\Vert^{#3}_{{#2}¥}}

\begin{document}
\title[Radon, cosine and sine transforms]
{Radon, cosine and sine transforms
on  Grassmannian manifolds
}
\author{Genkai Zhang
}
\address{Department of Mathematics, Chalmers University of Technology and
G\"oteborg University
 , G\"oteborg, Sweden}
\email{genkai@math.chalmers.se}

\thanks{Research  supported by the Swedish
Science Council (VR)\\
This is an updated version
of the article ArXiv 0810.5257 with the same title. I thank Siddhartha Sahi for pointing
out to me several misprints in the earlier version.
}

\def\bbK{\mathbb K}
\def\bKk{\mathbb K^k}
\def\bKn{\mathbb K^n}
\def\Cos{\text{Cos}}
\def\Sin{\text{Sin}}

\begin{abstract}
Let $G_{n,r}(\bbK)=G/K$ 
be the
 Grassmannian manifold
of $k$-dimensional $\bbK$-subspaces in $\bbK^n$
where $\bbK=\mathbb R, \mathbb C, \mathbb H$
is  the  field of real, complex or quaternionic numbers.
The Radon, 
cosine 
and sine transforms, 
$\mathcal R_{r^\prime, r}$, $\mathcal C_{r^\prime, r}$ and $\mathcal S_{r^\prime, r}$,
from  the
 $L^2$ space $L^2(G_{n,r}(\bbK))$ 
to the space  $L^2(G_{n,r^\prime}(\bbK))$,
for $r, r^\prime \le n-1$
are defined as integral operators in terms inclusion
relations of and angles between the subspaces.
 We compute  the  spectral symbols
of the transforms and characterize their images
under the decomposition
of $L^2$ spaces into irreducible subspaces of $G$.
 For that purpose
we prove two Bernstein-Sato type formulas
on general root systems  of type BC for
the sine and cosine type functions.
It is observed further that the Knapp-Stein intertwining
operator for certain induced representations
is given by the sine transform and as a consequence
we give
the unitary structure of the Stein's complementary
series  in the compact picture.

\end{abstract}

\maketitle








\def\mP{\mathcal P}
\def\mF{\mathcal F}
\def\SdN{\mathcal S^N}
\def\eSdN{\mathcal S_N^\prime}
\def\Gc{G_{\mathbb C}}
\def\rqs{G\backslash K_{\mathbb C}}

\baselineskip 1.25pc

\section{Introduction}

\newtheorem{theo+0}    {Theorem}    
\newenvironment{theo-abc}{\begin{theo+0}}{\end{theo+0}}

The present paper is a continuation
of our earlier papers \cite{gz-imrn}
and \cite{gz-radon-tams} on Radon and related transforms.
We shall give here a rather unified approach to the
 Radon, cosine and sine transforms 
on  Grassmannian manifolds.
Let $\mathbb K=\br, \bc, \mathbb H$ be the field of real, complex or quaternionic
numbers, and
$G_{n,r}=G_{n, r}(\mathbb K)$  be the Grassmannian manifold 
of $k$-dimensional subspaces over $\mathbb K$ in $\mathbb K^n$.
For $1\le r \le r^\prime \le n-1$, $\nu\ge 0$, the Radon, cosine and sine
transforms $\mathcal R=\mathcal R_{ r^\prime, r },
\mathcal C=\mathcal C_{ r^\prime, r }^{(2\nu)},
\mathcal S=\mathcal S_{ r^\prime, r }^{(2\nu)}: 
C^{\infty}(G_{n, r})
\to  C^{\infty}(G_{n, r^\prime}
)$ are defined,
for $\eta\in 
G_{n, r^\prime}$,
 by
\begin{equation}\label{rad-def}
(\mathcal Rf)
(\eta)=\int_{\xi\subset \eta}
f(\xi)
d_\eta \xi, 
\end{equation}
\begin{equation}\label{cos-1}
(\mathcal Cf)
(\eta)=\int_{G_{n, r}} |\Cos(\xi, \eta)|^{2\nu}
f(\xi)
d_\eta \xi, 
\end{equation}
\begin{equation}\label{sin-1}
(\mathcal Sf)
(\eta)=\int_{G_{n, r}}
|\Sin(\xi, \eta)|^{2\nu}
f(\xi)
d_\eta \xi, 
\end{equation}
where $d_\eta \xi$ is certain probability
measure on the set $\{\xi\in G_{n, r}: \xi\subset \eta\}$ 
invariant with
respect to the group of $\bbK$-unitary transformations
of $\eta$. The definition of the sine and cosine functions
is given in Definition 3.1. Now the Grassmannian
manifold $G_{n,r}$ is 
a compact Riemannian symmetric  space $G_{n, r}=G/K$, 
and the space $C^\infty(G_{n,r})$ is
decomposed under $G$ into irreducible
subspaces  by the Cartan-Helgason theorem.
The operators $\mathcal R$,
$\mathcal C$ and $\mathcal S$
are $G$-invariant and hence
they act on the irreducible subspaces
as  scalars, which can also be interpreted
as the spectral symbols of the transforms.
The spectral symbol
of the Radon transform
has been found earlier
in \cite{Grinberg-grs}
by integral computations using
explicit formulas for highest weight vectors.
In the present paper we shall   find the spectral symbols of the cosine
and sine transforms and 
a characterization 
of their images.
We  prove (see Theorem 5.4)

\begin{theo+0} The images of the transforms $\mathcal C$ and $\mathcal S$
are given by (5.1)-(5.3) in terms of the decomposition of $L^2(G_{n, r^\prime})$
under $G$.
\end{theo+0}

Our work generalizes the  recent results of \cite{Alesker-Berns}
on cosine transform where the case  $\mathbb K=\mathbb R$ and $2\nu=1$ is considered.
As application we find further the unitary 
structure of the Stein's complementary series
of the group $GL(2r, \mathbb K)$
realizing on the Grassmannian $G_{2r, r}$; see Section 6 for the notations.

\begin{theo+0} The Knapp-Stein intertwining operator
for the induced representations $\text{Ind}_P^{GL_{2r}}(t)$
and $\text{Ind}_P^{GL_{2r}}(-t)$ is given by the sine transform $\mathcal J_t=\mathcal S_{\nu}$. The unitary structure of the Stein's complementary
series for $0<t< \frac 12$ is given by $(f, g)_{\nu}=\frac 1{N_\nu^\prime} 
(\mathcal 
J_t f, g)$
where $(\mathcal 
J_t f, g)$  is explicitly given by (6.7).
\end{theo+0} 

To compute the spectral symbols
we shall use
the spherical
transform on root systems.
We observed in \cite{gz-radon-tams} (and it is known
in the rank one case)
that the squared Radon transform
$\mathcal R_{r^\prime, r}^\ast
\mathcal R_{r^\prime, r}$,
on a non-compact  symmetric matrix domain
has an integral kernel being
 a product of hyperbolic sine functions
whereas \cite{gz-imrn} 
the  Berezin transform has a product
of hyperbolic cosine functions, both
 in terms of the geodesic coordinates.
The spectral symbol of the Berezin transform
is then the Harish-Chandra spherical transform
of the hyperbolic cosine functions $\cosh^\nu t$.
 We find in \cite{gz-imrn}
this   transform  on a root system of type BC
by deriving
first certain Bernstein-Sato type formulas
for hyperbolic cosine functions $\cosh^\nu t$,  performing
spherical transform  and then by using the
observation that the functions $\cosh^{-\nu} t$, $\nu \to \infty$,
tend to
the Dirac delta function at $t=0$.
(This idea  has been used  before and might be due to
Berezin.)
Using the techniques developed there we derive
in this paper some corresponding Bernstein-Sato type formulas
for the sine and cosine functions on the compact torus
and we obtain the Cherednik-Opdam transform
for the sine and cosine transforms.

\begin{theo+0} \begin{enumerate}
\item
There exists a Weyl group
invariant polynomial $\mathcal M_{\delta}$
of the trigonometric Cherednik 
operators on the torus 
$\mathbb R^r/{2\pi Q^\vee}$, such that
the following Bernstein-Sato type formula hold
$$
\mathcal M_{\delta}
|\Cos (t)|^{\delta}
=\prod_{j=1}^r(\delta +a(j-1))
(\delta +\iota -1 +a(r-j))
|\Cos (t)|^{\delta-2},
$$
$$
\mathcal M_{\delta}
|\Sin (t)|^{\delta}
=\prod_{j=1}^r(\delta +a(j-1))
(\delta +\iota +2b-1 +a(r-j))
|\Sin (t)|^{\delta-2}.
$$
\item
The spherical transforms
$c_{\nu}(\m)=\widehat{|\Cos (t)|^{2\nu}}(\m)$
and
$s_{\nu}(\m)=\widehat{|\Sin (t)|^{2\nu}}(\m)$
are given by
\begin{equation*}
c_{\nu}(\m)
=N_\nu \prod_{j=1}^r\frac{(\nu+1+\frac a2(j-1)-\frac{m_j}2)_{\frac {m_j}2}}
{(\nu+1+\iota +b+ a(r-1)-\frac a2(j-1))_{\frac {m_j}2}}
\end{equation*}
and
\begin{equation*}
s_{\nu}(\m)
=N_\nu^\prime 
\prod_{j=1}^r\frac{(\nu+1+\frac a2(j-1)-\frac{m_j}2)_{\frac {m_j}2}}
{(\nu+1+\iota +b+ a(r-1)-\frac a2(j-1))_{\frac {m_j}2}} 
\end{equation*}
for root systems of Type D or C.
 \end{enumerate}
\end{theo+0} 
See Theorem 4.1 and 4.3.
Our methods can also be applied  Radon transform,
giving a different proof of some of the main results of Grinberg \cite{Grinberg-grs},
but we will not present the details here.

We make some general remarks on
the three transforms and also on the Berezin transform.
 The cosine and
sine transforms can be factorized in
 terms of the Radon transform; see Lemma 3.4.
The integral kernels 
of the two transforms
are related, roughly speaking,
by a rotation of $\frac \pi 2$. 
On the other hand,
the square $R_{r^\prime, r}^\ast R_{r^\prime, r}$
for $r^\prime >r$ is  actually
the cosine transform; see also \cite{gz-radon-tams}
for the case of non-compact symmetric matrix domains.
Now there is a fourth transform, the Berezin transform
which appears
 naturally in the branching
rule of holomorphic representations and quantization; see 
\cite{gkz-bere-rbsd}, \cite{gz-imrn} and references therein.
We prove also that the cosine transform
on  Grassmannians is  the compact analogue
of the  Berezin transform \cite{gkz-bere-rbsd}
on non-compact symmetric spaces. 
As another application we find the branching rule
of certain scalar holomorphic representations
under the group $G$; see Section 6.2.

The sine transform $\mathcal S_{r, r}$
is related to
the unitary representations of a  large non-compact group
$GL(2r, \mathbb K)$.
In \cite{Stein-slnc} Stein proved
that there is a family of unitary irreducible
representations of 
$GL(2r, \mathbb C)$
induced from non-unitary representations
of a maximal parabolic subgroup, also called
the Stein's complementary series.
The unitarity of those representation was
proved in \cite{Stein-slnc} by computing
the Fourier transform
of the kernel $|\det(x)|^{\alpha}$
of the  Knapp-Stein
intertwining operator; it gives then
the unitary structure in the so-called
{\it non-compact picture} of the induced representations.
Vogan \cite{Vogan-gln}
proved that the Stein's results can be extended
to  $GL(2r, \mathbb K)$ and proved
the corresponding existence of the
series. See also \cite{
Howe-Lee-jfa,Sahi-crelle,gkz-ma,Alesker-Berns}
for related study and applications.

After this paper was finished I was
informed by Professor Semyon Alesker 
that Theorem 5.4 (1) was partly proved
by him in an unpublished preprint \cite{Alesker-preprint} (under the
assumption that $\nu$ is not a negative half-integer)
using different methods. It is my pleasure to thank him
for his comments on the earlier versions of this paper. I thank also Professor B. Rubin
for informing me that a special case
of  Proposition 5.1 was also obtained in \cite{O-R-pre}.
The expert and careful comments
on an earlier version of this paper by the referee
are greatly acknowledged.

\section{Grassmannian manifolds $G_{n, r}$
and the irreducible decomposition of $L^2(G_{n, r})$}

\subsection{Grassmannian manifolds $\mathcal X=G_{n, r}$
as symmetric spaces}

Let $\bbK=\br, \bc, \mathbb H$ be the field of real, complex
or quaternionic numbers with the standard involution (conjugation)
$x\to \bar x$ and let
$a=\dim_\br \bbK =1, 2, 4$. Let $M_{n, m}:=M_{n, m}(\mathbb K)$
be the space of all $n\times m$-matrices, also viewed
as $\bbK$-linear transformations from $\bbK^m$ to
$\bbK^n$, where $\bbK$ acts on the right by scalar multiplication.
Denote $x^\ast =\bar x^T$, the conjugated transpose,
for  $x\in M_{n, m}$.
 Let
$$G:=U(n, \bbK)=\{g\in M_{n, n}; g^\ast g=I_n\}=O(n), \, U(n),\, Sp(n)
$$
be the orthogonal, unitary, and symplectic  groups accordingly.
 Denote, for $x\in M_{n, n}$, 
$\det_{\mathbb R}(x)$ the determinant of $x$ as a real linear
transformation on $\bbK^n=\mathbb R^{an}$. (For $\bbK=\mathbb C$
$ \det_{\mathbb R}(x)= |\det_{\mathbb C}(x)|^2$ and 
for $\bbK=\mathbb H$,  $ \det_{\mathbb R}(x)=\det_{\mathbb H}(x)^4$
where $\det_{\mathbb H}(x)$ is the so-called Dieudonn\'e determinant.)

Consider, for $r\le n$,  the Stiefel manifold $S_{n, r}$   of
all orthonormal $r$-frames in $\mathbb K^n$, 
and 
$$
\mathcal X=G_{n, r}
$$
 the Grassmannian manifold
of all $r$-dimensional subspaces over $\bbK$.
$S_{n, r}$ 
can be realized as the set
of all matrices  $x\in M_{n, r}$ such that
$x^\ast x=I$, namely
 $\bbK$-linear isometric transformations $x\in \bbK^r\to \mathbb K^n$.
Each
$x\in S_{n, r}$ defines uniquely a $r$-dimensional
subspace $\xi$ over $\bbK$,
$$\xi=\{x\}:=x\mathbb K^k\subset \mathbb K^n \in  G_{n, r}.
$$
Thus $\mathcal X$  is identified with the space of orbits in $S_{n, r}$ 
under the action of the unitary group $U(r, \bbK)$
on $\mathbb K^r$,
$$
\mathcal X=G_{n, r}=S_{n, r}/U(r, \bbK).
$$ 

Let 
$$
\xi_0=\mathbb K^r\oplus 0=\{x_0\}\in G_{n, r}, \qquad 
x_0=\begin{bmatrix} I_k\\
0
\end{bmatrix}\in S_{n, r}.
$$
We will fix $
\xi_0$ and $x_0$ as reference points of $\mathcal X$ and $S_{n, r}$,
The group $G$ acts on 
$S_{n, r}$ and  $\mathcal X$ 
by the defining action;
they are then realized 
as a homogeneous and respectively symmetric space,
\begin{equation}
  \label{eq:sym-rl-gr}
S_{n, r}=G/U(n-r, \bbK), \quad 
\mathcal X=G/K=U(n, \bbK)/U(r, \bbK)\times U(n-r, \bbK).
\end{equation}
Here $K=U(r, \bbK)\times U(n-r, \bbK)
$
is the isotropic subgroup of $\xi_0$
consisting of
block diagonal matrices of the form
$$
\begin{bmatrix} A&0\\
0&D
\end{bmatrix}, \quad A\in U(r, \bbK), \quad D\in  U(n-r, \bbK).
$$

Due to the $G$-isometry  between $G_{n, r}$
and $G_{n, n-r}$ we can and will assume in
 this paper, if nothing else is
stated, that 
\begin{equation}\label{2r-leq-n}
2r\le n.
\end{equation}

\subsection{Irreducible decomposition of $L^2(\mathcal X)$}
To describe the irreducible decomposition of $L^2(\mathcal X)$ under $G$
we let $\fg$ be the Lie algebra of $G$ and 
$\fg=\fk+\fp$  the Cartan decomposition of $\fg$, where
$\fk$ is the Lie algebra of $K$. The linear subspace
 $\fp$  consists of $n\times n$ matrices of the block form
$$
p_X=\begin{bmatrix}
 0&X\\
X^\ast&0\end{bmatrix}
, \quad X\in M_{r, n-r},
$$
which will be identified with $M_{r, n-r}$ via the mapping $X\mapsto p_X$. 
We let 
$\fa$ be the linear subspace of $\fp=M_{r, n-r}$ 
consisting of matrices of the form
$$
X=\begin{bmatrix}
\diag(t_1, \cdots, t_r)  & \, 0_{r, n-2r}\end{bmatrix}
=t_1E_1+\cdots +t_r E_r,  \quad t_1, \cdots, t_r\in \mathbb R
$$
with  $E_j$ being the matrix having $1$ on the $(j, j)$ position
and $0$ on the rest positions, $j=1, \cdots, r$. The maximal
torus $A=\exp(\fa)$ is then the group
 of $\exp(X)$, with
\begin{equation}
  \label{eq:expX}
\exp(X)=\begin{bmatrix}\diag(\cos t_1,\cdots, \cos t_r) &
\diag(\sin t_1,\cdots, \sin t_r) &0\\
\diag(-\sin t_1,\cdots, -\sin t_r) &
\diag(\cos t_1,\cdots, \cos t_r) &0\\
0& 0 & I
\end{bmatrix},
\end{equation}
written as  block matrix under the decomposition of $\mathbb K^n =
\mathbb K^r \oplus \mathbb K^r \oplus \mathbb K^{n-2r}$.
Let $\fg^\ast= \fk +i\fp=\fk +\sqrt{-1}\fp
$, the non-compact dual of $\fg$.
 The root system of $R(\fg^\ast, i\fa)$ is of the form
$$
R(\fg^\ast, i\fa)=
\{\pm \varepsi_j \pm \varepsi_k\}
\cup \{\pm \varepsi_j \} \cup
 \{\pm 2\varepsi_j \} 
$$
with respective multiplicities $a$, $a-1$ and $a(n-2r)$. Here $\{\varepsi_j\}$
is the dual basis of $\{iE_j\}$. It is of type C (if $n=2r$, $a>1$),
 type BC (if $n>2r$, $a>1$), 
 type B (if $n>2r$, $a=1$), or
 type D (if $n=2r$, $a=1$). It is understood here
and in Section 4 that
roots with zero multiplicity will not appear,
so that all types are considered as special case of type BC.
 We fix an ordering
so that  $\varepsi_1 >\cdots >\varepsi_r>0$
(and the condition $\varepsi_r >0$ is dropped for type D).

The following result follows from the Cartan-Helgason theorem
\cite{He2}. (Strictly speaking the theorem
is for $G/K$ with simply connected and semisimple Lie group $G$. 
The group $G=U(n, \mathbb K)$ for $K=\mathbb C$ or $\mathbb H$
is not semisimple, but the result can be reduced
to the semisimple group $SU(n, \mathbb K)$. However
for $\mathbb K=\mathbb R$ the group $G=O(n)$
is not connected nor $SO(n)$ is simple connected. This
imposes more conditions on the highest weights
than those given by that Theorem.)

\begin{prop+}  Under the action of $G=U(n, \mathbb K)$ the space $L^2(\mathcal X
)$ decomposes
as
   \begin{equation}
    \label{eq:L2-Gnr}
L^2(\mathcal X)=\sum_{\m}^{\oplus}     V^{\m}
  \end{equation}
with multiplicity free, where each $V^{\m}$ has highest weight $\m=m_1\varepsi_1
+\cdots m_r \varepsi_r$ with $\frac{m_j}2$ being integers and
$$
m_1\ge m_2\cdots \ge m_r\ge 0
$$
for $\mathbb K=\mathbb C, \mathbb H$
or $\mathbb K=\mathbb R$ with $n-2r>0$ and
$$
m_1\ge m_2\cdots \ge |m_r|
$$
for $\mathbb K=\mathbb R$
and $n=2r$. Each $V^{\m}$ has a unique
$K$-invariant function $\phi_{\m}$, called
spherical polynomial,  normalized so that
$\phi_{\m}(\xi_0)=1$.
\end{prop+}

We will view the spherical polynomial
$\phi_{\m}$ on $\mathcal X$ as defined on $\fa=\mathbb R^r$ via the
exponential mapping $\fa\to G_{n, r}$,
$X\to \exp(X)\cdot \xi_0$, namely
$$
\phi_{\m} (t_1, \cdots, t_r):=\phi_{\m}
(\exp(t_1E_1 + \cdots  +t_rE_r)\cdot \xi_0).
$$

Finally we recall the polar decomposition of $G_{n, r}$. 
For any function $f$ on $G_{n,r}$ we have
the following formula (which in turn fixes
a normalization of the invariant measure)
\begin{equation}
  \label{eq:polar}
\int_{G_{n, r}} f(\xi)d\xi =\int_K \int_{A} f(k a)\cdot \xi_0)dkd\mu(a)  
\end{equation}
where
$$
d\mu(a)
=d\mu_{R}(t)=\prod_{\a\in R_+} |2\sin \alpha(it))|^{m_{\a}}
dt_1\cdots dt_r, \quad a=\exp(t)=\exp(t_1 E_1 +\cdots +t_r E_r)
$$
with $m_{\a}$ the root multiplicity of $\a$ 
and $dk$  the normalized measure on $K$;
see \cite[Chapter I, Theorem 5.10]{He2}.
\section{Radon,  cosine and sine transforms}

\subsection{Cosine and sine of angles between
subspaces}
We fix the standard Euclidean norm  on $\mathbb K^n=\mathbb R^{an}$. For any
convex subset $S$ in a
$d$-dimensional real subspace of $\mathbb K^n$ we let $\vol_{d}(S)$ be the corresponding Euclidean volume.  
The following definition is given in \cite{Alesker-Berns}
for the real Grassmannians.
\begin{defi+} 
\label{cos-sin-def} 
If $r\le r^\prime$ the cosine of the angle between two subspaces
$\xi\in G_{n, r}$ and  $\eta\in G_{n, r^\prime}$
is defined by
\begin{equation}
  \label{eq:cos-def}
|\Cos(\xi, \eta)|=
\left(\frac{\vol_{ra}P_{\eta }(A) }{\vol_{ra}(A)}\right)^{\frac 1a},
\end{equation}
where $P_\eta$ is the orthogonal projection
from $\mathbb K^n $ onto $\eta\subset \mathbb K^n$
and $A\subset \xi$ is any convex set of non-zero volume.
If $r\le n-r^\prime$ the sine of the angle 
 between two planes 
$\xi\in G_{n, r}$ and  $\eta\in G_{n, r^\prime}$ is defined by
$$
|\Sin(\xi, \eta)|=|\Cos(\xi, \eta^{\perp})|.
$$
For general $r$ and  $r^\prime$ we
define $|\Cos(\xi, \eta)|$ and $|\Sin(\xi, \eta)|$ by the symmetry
condition
$$
|\Cos(\xi, \eta)| =|\Cos(\eta, \xi)|,\qquad 
 |\Sin(\xi, \eta)|= |\Sin(\eta, \xi)|.
$$
\end{defi+}

 In \cite{Grinberg-Rubin} Grinberg and Rubin introduce
also a cosine of angle,  
$\text{COS}^2(y, x)$, between two elements 
$x\in S_{n, r}$ 
and $y\in S_{n, r^\prime}$, for $r\le r^\prime$, defined to be the $r\times r$-semi-positive
definite matrix
\begin{equation}
  \label{eq:GR-cos}
\text{COS}^2(y, x)=x^\ast y y^\ast x.  
\end{equation}

The following lemma computes the cosine (3.1)
in terms of   (\ref{eq:GR-cos})
\begin{lemm+}
\begin{enumerate}
\item Let $1\le r\le r^\prime$ and the notations be as above.
The two
 cosine functions
 (\ref{eq:cos-def}) and    (\ref{eq:GR-cos})
are related by,
$$
|\Cos (\eta, \xi)|^{a}=(\det_{\mathbb R}\text{COS}^2(y, x)  )^{\frac 12}
$$
where $\eta=y\mathbb K^{r^\prime}\in G_{n, r^\prime}$,
$\xi=y\mathbb K^{r}\in G_{n, r}$. 
\item Let $1\le r = r^\prime \le n-r$. 
Write  $\xi\in G_{n, r}$ as $\xi =k\exp(
t_1E_1+\cdots +t_r E_r)\cdot \xi_0$, $k\in K$. Then
  \begin{equation}
|\Cos (\xi, \xi_0)|^{\delta} =\prod_{j=1}^r |\cos t_j|^\delta.
  \end{equation}
  \begin{equation}
|\Sin (\xi, \xi_0)|^{\delta} =\prod_{j=1}^r |\sin t_j|^\delta.
  \end{equation}
\end{enumerate}
\end{lemm+}
\begin{proof} We recall first that if $T: \mathbb R^p \to 
\mathbb R^q$, $p\le q$, is a linear transformation and $B\subset \mathbb R^p$ is
any convex set of positive volume, then $\vol_{p}(T(B))=\det(T^t T)^{\frac 12}
\vol_p (B)$, where $T^t$ is the transpose of $T$ with respect
to the Euclidean  metrics in $\mathbb R^p$ and $\mathbb R^q$.
 Indeed, let $T=U(T^t T)^{\frac 12}$
 be the polar decomposition
of $T$ with $U$ being a partial
isometry. If $\text{rank}(T)< p$, then both
$\vol_{p}(T(B))$ and $\det(T^t T)^{\frac 12}$ are zero and the formula
is trivially true.  If $\text{rank}(T)= p$ then $U^tU=I$  and $U$
is an isometry. We have then
$$
\vol_p(T(B))=\vol_p (U(T^t T)^{\frac 12}(B))
=\vol_p ((T^t T)^{\frac 12}(B))
=\det ((T^t T)^{\frac 12})\, \vol_p(B),
$$
proving the identity. 

Now if $\xi=\{x\}=x\mathbb K^r$ and $\eta=\{y\}=y\mathbb K^{r^\prime}$,
with isometries $x\in S_{n, r}$ and $y\in S_{n, r^\prime}$,
we have $P_\eta=yy^\ast$. Let $A\subset \xi=
x\mathbb K^r$
 be
any convex set of nonzero volume. We write $A 
=x (B)$ with $B\subset \mathbb K^r$ a convex set,
and $\vol_{ra}(B)=\vol_{ra}(A)$ since $x$ is an isometry.
Its image under $P_\eta$
is $P_\eta(A)
=yy^\ast x(B)$. We apply the previous formula with $T =P_\eta x=yy^\ast x$,
noticing that
$
T^\ast T=x^\ast P_\eta^2 x =x^\ast P_\eta x =x^\ast  yy^\ast x,
$
$$
\vol_{ra}P_{\eta }(A) 
=\left(\det_{\mathbb R}(x^\ast  yy^\ast x)\right)^{\frac 12}
\vol_{ra}(A)
=\det_{\mathbb R}(\text{COS}^2(x, y))^{\frac 12}
\vol_{ra}(A).
$$
This proves the first part, and the second part
is  then a special case
by using the formula (\ref{eq:expX}).
\end{proof}

\subsection{Factorization and diagonalization
of the  cosine and sine transforms}
We define the Radon transform
$\mathcal R_{r^\prime, r}:
 C^\infty(G_{n, r})\to C^\infty(G_{n, r^\prime})$ by
\begin{equation}\label{def-radon-1}
(\mathcal R_{r^\prime, r}f)(\eta)=\int_{\xi\in G_{n, r}; \xi\subset\eta} 
f(\xi)d_\eta\xi,
\end{equation}
if $r < r^\prime$ and 
\begin{equation}
\label{def-radon-2}
(\mathcal R_{r^\prime, r}f)(\eta)=
\mathcal R_{r^\prime, r}^\ast  f(\eta)=
\int_{\xi\in G_{n, r}; \xi\supset \eta} f(\xi) d_\eta\xi. 
\end{equation}
if $r > r^\prime$. Here $d_\eta \xi$ is the unique Riemannian measure
on the subset induced from the Riemannian measure  on $G_{n, r}$;
see \cite{Grinberg-Rubin} for an expression of
$\mathcal R_{r^\prime, r}$ and the measure in terms of
the  realization (\ref{eq:sym-rl-gr}). In particular
the Radon transform $\mathcal R_{r^\prime, r}$
commutes with the actions of $G$ on
$ C^\infty(G_{n, r})$ and  $C^\infty(G_{n, r^\prime})$.

\begin{defi+}  Let  $1\le r, r^\prime< n$ and $\nu \ge 0$. We define the
sine and cosine transforms from $C^\infty(G_{n, r})$
to $C(G_{n, r^\prime})$, by
$$
\mathcal S^{(\nu)}_{r^\prime, r}f(\eta)=\int_{G_{n, r}}
|\Sin (\eta, \xi)|^{2\nu} f(\xi)  d\xi,
\quad
\mathcal C^{(2\nu)}_{r^\prime, r}f(\eta)=\int_{G_{n, r}}
|\Cos (\eta, \xi)|^{2\nu} f(\xi)  d\xi.
$$
\end{defi+}

The sine and cosine transforms are related to the Radon transform
via the following Lemma.
This Lemma in the case when $2\nu=1$ 
and  $\mathbb K=\mathbb R$ is 
known
and is proved in Lemma 1.7 in \cite{Alesker-Berns}.
The same method can be applied to the
present case for general $\nu \ge 0$,
   by using a variant  of the Cauchy-Kubota
formula 
$$
\vol_{ar}(B)^\delta= c_\delta \int_{\xi\in G_{r^\prime, r}} 
\vol_{ar}(\Pr_{\xi} B)^\delta d\xi, 
$$
for any convex polygon in $B\subset \xi_0$, which
follows easily by the invariance of both sides
under translations and linear actions on $\xi_0$.
We skip the elementary proof.

\begin{lemm+}
 Let $1\le r\le r^\prime < n$, and
$\nu\ge 0$. The cosine and sine transforms
can be factorized as
$$
\mathcal C^{(\nu)}_{r, r^\prime}
=c_1
\mathcal C^{(\nu)}_{r, r} \mathcal R_{r, r^\prime}, \qquad
\mathcal S^{(\nu)}_{r, r^\prime}=c_2
\mathcal S^{(\nu)}_{r, r} \mathcal R_{r, r^\prime},
$$
where $c_1=c_1(\nu)$ and $c_2=c_2(\nu)$ are some positive 
constants.
 \end{lemm+}
The constant $c_1$ and $c_2$ can be computed by using
the integral formula in \cite{gz-radon} but we will not need it here.

By taking conjugate of the above formulas we 
get  factorizations
 for any $1\le r, r^\prime <n$, noticing
that $(\mathcal C_{r, r^\prime}^{(\nu)})^\ast
=\mathcal C_{r^\prime, r}^{(\nu)}$, 
 $(\mathcal S_{r, r^\prime}^{(\nu)})^\ast
=\mathcal S_{r^\prime, r}^{(\nu)}$. 

The transforms $\mathcal C^{(\nu)}_{r, r^\prime}$
and  $\mathcal S^{(\nu)}_{r, r^\prime}$ clearly intertwine the action of
$G$.
Consider the corresponding decomposition of $L^2(G_{n, r^\prime})$
according to Proposition 2.1,
  \begin{equation}
    \label{eq:L2-Gnr-pr}
L^2(G_{n, r^\prime})=\sum_{\m}^{\oplus} W^{\m}.
  \end{equation} 
Thus $\mathcal C^{(\nu)}_{r, r^\prime}$
and  $\mathcal S^{(\nu)}_{r, r^\prime}$ are diagonal operators
up to a normalization.
The eigenvalue of 
$\mathcal R_{r, r^\prime}$ has been found earlier by Grinberg \cite{Grinberg-grs};
see Theorem 5.2 below.
In view of Lemma 3.4 above,  to find the eigenvalues of
 $\mathcal C^{(2\nu)}_{r, r^\prime}$ and  $\mathcal S^{(2\nu)}_{r, r^\prime}$
we need only consider the case when $r^\prime =r$.

Denote $c_{\nu}(\m)=
c_{\nu, r}(\m)=
$ and
$s_{\nu}(\m)$ 
the eigenvalue of 
$\mathcal C^{(\nu)}_{r, r}$ 
and respectively $\mathcal S^{(\nu)}_{r, r}$ 
on $V^{\m}$. They can be  evaluated by 
  \begin{equation}
\mathcal C^{(\nu)}_{r, r}
\phi_{\m}(\xi_0)=c_{\nu}(\m)
\phi_{\m}(\xi_0)=c_{\nu}(\m).
  \end{equation} 
Namely
  \begin{equation}
c_{\nu}(\m)=
\widehat{|\Cos (\cdot, \xi_0)|^{2\nu}}(\m)
=\int_{G_{n, r}}|\Cos (\xi, \xi_0)|^{2\nu} \phi_{\m}(\xi)d\xi,
  \end{equation} 
is the spherical transform of 
$|\Cos (\cdot, \xi_0)|^{2\nu}$.
Similarly for 
$s_{\nu}(\m)$.
They are  integrals of  $K$-invariant functions. We use
the polar coordinates (\ref{eq:polar}),   which
 further can be written as an integral
on the quotient  $\mathbb R^r/{2\pi Q^\vee}$
of $\mathbb R^r$  by the 
(spherical) coroots lattice $Q^\vee$, 
namely the lattice generated by  $\alpha^\vee
=\frac{\hat \alpha}{(\alpha, \alpha)}$
for $\alpha\in  R$, where $\hat \alpha$ is
the dual of $\alpha$, 
$\lam(\hat \alpha)=(\lam, \alpha)$; see \cite[
Chapter I, Theorem 5.10]{He2}.

\begin{lemm+} 
Denote
$$
|\Cos (t)|
:=|\prod_{j=1}^r \cos t_j|, \qquad
|\Sin (t)|
:=|\prod_{j=1}^r \sin t_j|,
$$
with 
$$
\quad t=t_1E_1 + \cdots+ t_rE_r=
(t_1, \cdots, t_r)\in \mathbb R^r/{2\pi Q^\vee}.
$$
Let $\nu\ge 0$. The eigenvalues
$c_{\nu}(\m)$ and 
$s_{\nu}(\m)$ of the cosine respectively sine
transforms are given by
the spherical transforms
  \begin{equation}
c_{\nu}(\m)=
\widehat{|\Cos (t)|^{2\nu}}(\m)
=\int_{\mathbb R^r/{2\pi Q^\vee}}
|\Cos (t)|^{2\nu}
 \phi_{\m}
(t_1, \cdots, t_r )
d\mu_R(t),
  \end{equation} 
  \begin{equation}
s_{\nu}(\m)
=\widehat{|\Sin (t)|^{2\nu}}(\m)
=\int_{\mathbb R^r/{2\pi Q^\vee}
}
|\Sin (t)|^{2\nu}
\phi_{\m}
(t_1, \cdots, t_r )
d\mu_R(t).
  \end{equation} 
\end{lemm+}

\section{Spherical transform on compact
torus associated with general root system of Type BC}

\subsection{Bernstein-Sato type formulas for cosine and sine functions
}
We consider the root system $R$ on a Euclidean space $i\fa=i\mathbb R^r$
($i=\sqrt{-1}$)
of type BC
\begin{equation}
\label{rt-sys}
R=\{\pm\varepsi_j \pm \varepsi_k, j\ne k\} 
\cup\{\pm \varepsi_k\}
\cup\{\pm 2\varepsi_k\}
\end{equation}
with general non-negative multiplicities $a$, $2b$ and $\iota$ for
the respective sets of roots. Here $\{\varepsi_j\}$
is the dual basis in $(i\fa)^\ast$ of
a fixed orthonormal basis $\{iE_j\}$ of 
 $i\fa$; the notation $\{E_j\}$
here coincides with that in Section 2.2.
We use the same 
convention there, so that type D and type B
will be considered as  special cases.

Denote by $W$
the Weyl group. We will compute the
spherical  transform of certain Weyl group invariant
sine and cosine
functions, by using the harmonic analysis of 
the Cherednik operators developed by Opdam \cite{Opdam-acta}. 
We follow the presentation there loc. cit., however
with our roots being twice of the roots there  and our
multiplicities  half of the ones there.

We let $D_j=D_{iE_j}$, $j=1\cdots, r$,
be the trigonometric Cherednik operators acting
on  functions on $\mathbb R^r/{2\pi Q^\vee}$,
\begin{equation*}
\begin{split}
D_j&=\partial_j -ia\sum_{k<j}\frac{1}{1-e^{-2i(t_k-t_j)}}(1-s_{kj})
+i a\sum_{j<k}\frac{1}{1-e^{-2i(t_j-t_k)}}(1-s_{jk})+
\\
&+ia \sum_{k\ne j}\frac{1}{1-e^{-2 i (t_j+t_k)}}(1-\sig_{jk})
+2i\iota \frac{1}{1-e^{-4 i t_j}}(1-\sig_{j})  \\
&+2ib \frac{1}
{1-e^{-2 i t_j}}(1-\sig_{j})
-i\rho_j,
\end{split}
\end{equation*}
where $s_{kj}$, $\sig_{kj}$ and $\sig_j$
are the reflections corresponding to the
roots $\varepsi_{j}-\varepsi_{k}$,
$\varepsi_{j}+\varepsi_{k}$,and
$\varepsi_{j}$. Here 
$$
\rho=\frac 12 \sum_{\alpha\in R^+} m_{\alpha} \alpha
=\sum_{j=1}^r \rho_j \varepsi_j, \quad
\rho_j=\iota+ b + a(r-j).
$$
is the half sum of positive roots, $m_{\alpha}$ being
the root multiplicities.
Let $\phi_{\m}$ be
the Heckman-Opdam Jacobi polynomials on the root
system (\cite{Heckman-Opdam-1}, \cite{Heckman-Opdam-2}, and \cite{Opdam-acta}),
 and 
$$
\hat f(\m):=\int_{\a/{2\pi Q^\vee}}
f(t)\phi_{\m}(t) d\mu_{R}(t)
$$
be the spherical (or  Jacobi) transform. 
Here  $Q^\vee$ is as before  the (spherical) coroot lattice.
Our objective
is to find the spherical  transform of the  functions
$|\Cos (t)|^{2\nu}$ and $|\Sin (t)|^{2\nu}$ in Lemma 3.5. 
We establish first certain Bernstein-Sato type formulas,
more exactly we will find certain Weyl group invariant
polynomials of the Cherednik operators mapping
$|\Cos t|^\delta$ to $|\Cos t|^{\delta-2}$.

\begin{theo+} \label{thm-1} 
Let $\delta\ge 0$ and $\mathcal M_{\delta}$ be the  operator
$$
\mathcal M_{\delta}:
\mathcal =
\prod_{j=1}^r (D_j^2 +(\delta+\rho_1)^2).
$$
Then the following Bernstein-Sato type formulas hold,
$$
\mathcal M_{\delta}
|\Cos (t)|^{\delta}
=\prod_{j=1}^r(\delta +a(j-1))
(\delta +\iota -1 +a(r-j))
|\Cos (t)|^{\delta-2}
$$
and
$$
\mathcal M_{\delta}
|\Sin (t)|^{\delta}
=\prod_{j=1}^r(\delta +a(j-1))
(\delta +\iota +2b-1 +a(r-j))
|\Sin (t)|^{\delta-2}.
$$
 \end{theo+}

\begin{proof} In \cite[Theorem 2.1]{gz-imrn} 
and \cite[Theorem 3.1]{gz-radon-tams} the following
formulas are proved for the 
hyperbolic sine and cosine functions, defined on $i\mathfrak a$,
\begin{equation*}
\begin{split}
&\quad \prod_{j=1}^r\big((iD_j)^2-(\delta+\rho(\xi_1))^2\big)
(\prod_{j=1}^r \cosh x_j)^{\delta}\\
&=\prod_{j=1}^r \big(\delta +a(j-1) \big)
\big( 1-\delta-\iota-a(r-j) \big)
(\prod_{j=1}^r \cosh x_j)^{\delta-2},
\quad x=(x_1, \cdots, x_r)\in \mathbb R^r,
\end{split}
\end{equation*}
\begin{equation*}
\begin{split}
&\quad \prod_{j=1}^r
\big((iD_j)^2-(\delta+\rho_1)^2\big)
(\prod_{j=1}^r \sinh x_j)^{\delta}\\
&=\prod_{j=1}^r \big(\delta +a(j-1) \big)
\big( \delta-1 +\iota +2b +a(r-j)\big)
(\prod_{j=1}^r \sinh x_j)^{\delta-2}, \quad x=(x_1, \cdots, x_r)\in \mathbb R_{+}^r,
\end{split}
\end{equation*}
where $iD_j$ is acting on the variable $x$.
The function
$(\sinh x)^{\delta}$ and $(\cosh x)^\delta$ are analytic
functions on the right
hand plane $\Re x>0$, 
so are  the
products  $\prod_{j=1}^r
(\sinh x_j)^{\delta}$
and $\prod_{j=1}^r
(\cosh x_j)^{\delta}$ on the product
of the right half plane 
 $\{x; \Re x_j>0\}$. All the identities
has a limit at the  points $x=it$, $t\in \mathbb R^r$, $t_j\ne 0, \frac{\pi}{2}$.
Our result follows then by taking the limit, observing
also that $\sinh^2(it_j)=(-1)\sin^2 t$.
\end{proof}

\begin{rema+} The above theorem can also be
proved directly by a straightforward but long
computations. First, we can choose a dense open
subset of the orbit $W\backslash\mathbb R^r/{2\pi Q^\vee}$
of the Weyl group so that the functions $\sin t_j >0, \cos t_j>0$.
Furthermore the operator $\mathcal M_{\delta}$
can be factorized as follows,
\begin{equation}
\label{factor}
\mathcal M_{\delta}=
\prod_{j=1}^r (D_j^2 +(\delta+\rho_1)^2)=
\prod_{j=1}^r (-iD_j +(\delta+\rho_1))
(iD_j +(\delta+\rho_1))
\end{equation}
and we may compute successively the action of
each factors.
 We have, for each fixed $j$,
\begin{equation*}
\prod_{l=1}^j(iD_l +(\delta+\rho_1)) \Cos (t)^{\delta}
=\prod_{l=1}^j (\delta +a(l-1))
\left(\Cos (t)^{\delta}
\prod_{l=1}^j\frac{e^{-it_l}}{\cos t_l}\right),
\end{equation*}
\begin{equation*}
\begin{split}
&\quad
 \prod_{l=j}^r (-iD_l +(\delta+\rho_1))
\left(\Cos (t)^{\delta}\prod_{l=1}^j\frac{e^{-it_l}}{\cos t_l}\right)\\
&
=\prod_{l=j}^r (\delta-1 +\iota +a(r-l))
\left(\Cos (t)^{\delta}\prod_{l=1}^j\frac{e^{-it_l}}{\cos t_l}\right)
\prod_{l=j}^r \frac {e^{it_l}}{\cos t_l}
\end{split}
\end{equation*}
\begin{equation*}
\prod_{l=1}^j(iD_l +(\delta+\rho_1))\Sin (t)^{\delta}
=i^j\prod_{l=1}^j (\delta +a(l-1))
\left(
\Sin (t)^{\delta}
\prod_{l=1}^j\frac {e^{-it_l}}{\sin t_l}
\right)
,
\end{equation*}
\begin{equation*}
\begin{split}
&\quad
\prod_{l=j}^r (-iD_l +(\delta+\rho_1))
\left(
\Sin (t)^{\delta}
\prod_{l=1}^j\frac {e^{-it_l}}{\sin t_l}
\right)\\
&
=(-i)^{r-j+1}
\prod_{l=j}^r (\delta-1 +\iota +2b +a(r-l))
\left(
\Sin (t)^{\delta}
\prod_{l=1}^j\frac {e^{-it_l}}{\sin t_l}
\right)\prod_{l=j}^r \frac {e^{it_l}}{\sin t_l},
\end{split}
\end{equation*}
by similar  computations as in \cite{gz-imrn} 
and \cite{gz-radon-tams}, which together with the factorization
(\ref{factor}) implies our theorem. (The pattern of
the identities is roughly that
each action of 
$( iD_j +(\delta+\rho_1))$ on the functions $\Cos \,t$
produces an extra factor $\frac{e^{it_j}}{\cos t_j}$,
further action by
$( iD_j +(\delta+\rho_1))$ produces
an extra factor $\frac{e^{-it_j}}{\cos t_j}$;
the $e^{\pm i t_j}$ cancels and we get a single factor of $\Cos^{-2}t$,
namely our formula.)
The above family of identities
is to be compared with the trivial trigonometric formulas
\begin{equation*}
(i\frac{d}{dt} +\delta) \cos (t)^{\delta}
=\delta \cos (t)^{\delta}
\frac{e^{-it}}{\cos t},
\end{equation*}
\begin{equation*}
(-i\frac{d}{dt} +\delta) 
\left(\cos (t)^{\delta}
\frac{e^{-it}}{\cos t}\right)
=(\delta-1)
\left(
\cos (t)^{\delta}
\frac{e^{-it}}{\cos t}\right)
\frac{e^{it}}{\cos t},
\end{equation*}
\begin{equation*}
(i\frac{d}{dt} +\delta) \sin (t)^{\delta}
=i\delta 
\left(
\sin (t)^{\delta}
\frac{e^{-it}}{\sin t}
\right),
\end{equation*}
\begin{equation*}
(-i\frac{d}{dt} +\delta)
\left(
\sin (t)^{\delta}
\frac{e^{-it}}{\sin t}
\right)
=(-i)(\delta -1)
\left(
\sin (t)^{\delta}
\frac{e^{-it}}{\sin t}
\right)\frac{e^{it}}{\sin t}.
\end{equation*}
It would be interesting to reformulate the identities
systematically in terms of Hecke algebras \cite{Cherednik-imrn97}.
\end{rema+}

\subsection{Spherical  transform for cosine and sine functions
}
We let $ N_\nu$ and $ N_\nu^\prime$ be the following normalization
constants,
\begin{equation}
  \label{eq:n-nu}
  N_\nu=\int_{\mathbb R^r
/{2\pi Q^{\vee}}
}
|\Cos|^{2\nu}(t) d\mu(t)
\end{equation}
and 
\begin{equation}
  \label{eq:n-nu-'}
  N_\nu^\prime
=\int_{\mathbb R^r/{2\pi Q^{\vee}}}
|\Sin|^{2\nu}(t) d\mu(t).
\end{equation}
Their exact values can be evaluated by using the Macdonald
formula for generalized Beta-integrals
(see \cite[Ex. 7, Sect. 10, Chapt. VII]{Macd-book}), 
\begin{equation}
\begin{split}
N_\nu& =2^{ar(r-1) +2rb +2r\iota}r!\prod_{1\le i<j\le r}
\frac{\Gamma(\frac a2(j-i +1)}
{\Gamma(\frac a 2(j-i)}\times\\
&\quad
\frac{
\Gamma_a( 1+b +\frac{\iota -1}2 +\frac a2(r-1)  )
\Gamma_a( \nu +1+\frac{\iota -1}2 + \frac a2 (r-1))
}
{
\Gamma_a( \nu +1+b+\iota + a(r-1) )
}
\end{split}
\end{equation}
\begin{equation}
\begin{split}
N_\nu^\prime &=2^{ar(r-1) +2rb +2r\iota}r!\prod_{1\le i<j\le r}
\frac{\Gamma(\frac a2(j-i +1)}
{\Gamma(\frac a 2(j-i)} \times
\\
&\quad
\frac{
\Gamma_a( 1 +\frac{\iota -1}2 +\frac a2(r-1)
 )
\Gamma_a( \nu +1+\frac{\iota -1}2 + \frac a2 (r-1)
)
}
{
\Gamma( \nu +1+b+\iota + a(r-1) -\frac a2(j-1) )
}.
\end{split}
\end{equation}
Here $\Gamma_a(\alpha)$ is the Gindikin's Gamma function
$$
\Gamma_a(\alpha)=\prod_{j=1}^r\Gamma(\alpha -\frac a2(j-1)).
$$
We recall also the Pochammer symbol $(\nu)_k =(\nu)(\nu +1)\cdots (\nu +k-1)$.
\begin{theo+}
 \label{thm-2}
The spherical  transforms
of the functions
$|\Cos (t)|^{2\nu}$
and $|\Sin (t)|^{2\nu}$
are given by
\begin{equation*}
c_{\nu, r}(\m):=
\widehat{|\Cos|^{2\nu}}
(\m)
=N_\nu \prod_{j=1}^r\frac{(\nu+1+\frac a2(j-1)-\frac{m_j}2)_{\frac {m_j}2}}
{(\nu+1+\iota +b+ a(r-1)-\frac a2(j-1))_{\frac {m_j}2}}
\end{equation*}
and
\begin{equation*}
s_{\nu, r}(\m):=\widehat{|\Sin|^{2\nu}}(\m)
=N_\nu^\prime 
\prod_{j=1}^r\frac{(\nu+1+\frac a2(j-1)-\frac{m_j}2)_{\frac {m_j}2}}
{(\nu+1+\iota +b+ a(r-1)-\frac a2(j-1))_{\frac {m_j}2}} 
\phi_{\m}(\frac \pi 2, \cdots, \frac \pi 2).
\end{equation*}
For root systems of Type C or Type D we have $b=0$,
$\phi_{\m}(\frac \pi 2, \cdots, \frac \pi 2)=\prod_{j=1}^r
 (-1)^{m_j}$ (see below),
and thus
\begin{equation*}
s_{\nu, r}(\m):=\widehat{|\Sin|^{2\nu}}(\m)
=N_\nu^\prime 
\prod_{j=1}^r\frac{(-\nu-\frac a2(j-1))_{\frac {m_j}2}}
{(\nu+1+\iota +b+ a(r-1)-\frac a2(j-1))_{\frac {m_j}2}}.
\end{equation*}
\end{theo+}

We need   the following elementary result, which
states simply that the normalized integration
of the cosine functions $\cos^n t$ 
and $\sin^n t$ tends
to the  $\delta$-function
at $t=0$ and
$t=\frac {\pi}2$ respectively.

 \begin{lemm+}
Suppose  $\phi$ be a bounded and continuous function
on $\mathbb R^r/{2\pi Q^\vee}$. Then
$$
\lim_{\nu\to \infty}\frac{1}
{N_{\nu}}
\int_{\mathbb R^r/{2\pi Q^\vee}
}
|\Cos|^{2\nu}(s) \phi(s) d\mu(s)
=\phi(0)
$$
and
$$
\lim_{\nu\to \infty}\frac{1}
{N_{\nu}^\prime}
\int_
{\mathbb R^r/{2\pi Q^\vee}} |\Sin|^{2\nu}(s)d\mu(s)
=\phi(\frac \pi 2, \cdots,
\frac \pi 2).
$$
 \end{lemm+}

The following lemma asserts that
$\phi_{\m}(\frac \pi 2, \cdots,\frac \pi 2)$ is always nonzero,
which is needed in Theorem 5.4.
 \begin{lemm+} 
 \begin{enumerate} 
\item
Suppose
$R$ is a root system of type  C, or D with general non-negative
 root
multiplicity. Then 
$$
\phi_{\m}(\frac \pi 2, \cdots,\frac \pi 2)
=\prod_{j=1}^r(-1)^{m_j}
$$
\item Suppose $R$ is the root system of the Grassmannian manifold
$G_{n, r}$. Then 
$$
\phi_{\m}(\frac \pi 2, \cdots,\frac \pi 2)
\ne 0.
$$
  \end{enumerate}  
 \end{lemm+}

 \begin{proof} (1) From the formula for $D_j$
we see that, if $R$ is a root system of type C, or D,
then 
$\{D_j\}$ are invariant under the map $(t_1, \cdots, t_r)
\mapsto (t_1+\frac\pi 2, \cdots, t_r +\frac\pi 2)
$.
Moreover it is easy to prove that 
the polynomial $f(t_1, \cdot, t_r):=
\phi_{\m}(t_1+\frac \pi 2, \cdots, t_r +\frac \pi 2)$
is also invariant under the Weyl group. We prove this
for type C, the  type D is  the same. 
The Weyl group is generated
by  the simple reflections $\sig_r$ 
and $s_{j, j+1}$, so we need only
to check the invariance for those elements. 
Notice that by definition
$\phi_{\m}$ is invariant under
the mapping $(t_1, \cdots, t_r)\mapsto (t_1, \cdots, t_r+\pi)$.
We have
 \begin{equation*}
 \begin{split}
(\sig_r f)(t_1, \cdots, t_r)&=
f(t_1, \cdots, -t_r)=\phi_{\m}(t_1+\frac \pi 2, \cdots, -t_r +\frac \pi 2)\\
&=\phi_{\m}(t_1+\frac \pi 2, \cdots, -(t_r +\frac \pi 2) +\pi)\\
&=\phi_{\m}(t_1+\frac \pi 2, \cdots, -(t_r +\frac \pi 2))\\
&=\phi_{\m}(t_1+\frac \pi 2, \cdots, (t_r +\frac \pi 2))\\
&=f(t_1, \cdots, t_r)
 \end{split}
 \end{equation*}
and  that $s_{j j+1} f=f$ is trivially true.
Thus $f$ is an eigenfunction of the Weyl group invariant
polynomials of the operators $\{D_j\}$ with
 the eigenvalues being the same 
as that of $\phi_{\m}$. Namely 
$f(t)=c\phi_{\m}(t)$ for some constant $c$, by the uniqueness
of the spherical eigenfunctions 
\cite{Opdam-acta}. Comparing the leading coefficients we
find the constant $c$, and the evaluation of $f$ at $t=0$
proves our result.

(2) Consider the case of  Grassmannian manifolds. We
need only to treat the case when
the corresponding   root system $R$ is of type B or BC,
namely $G_{n, r}(\mathbb K)$ for $n-r>r$. Write
$$
\xi_1:=\exp(\frac{\pi}2 E)
\cdot \xi_0=0\oplus  \mathbb K^{r}\oplus 0 \in G_{n, r},
$$
where
$\exp(\frac{\pi}2 E)$ is the short-hand notation
$\exp(\frac{\pi}2 E)
=\exp(\frac{\pi}2E_1 +\cdots +
\frac{\pi}2 E_r)$.
We have, by the integral formula
for spherical polynomials, \cite[Chapter IV, Proposition 2.2]{He2}, 
that
\begin{equation}
  \label{eq:int-form-sp}
\phi_{\m}(\xi_1)^2
=\int_{K}
\phi_{\m}(\exp(\frac{\pi}2 E)
k
\exp(\frac{\pi}2 E)
\cdot \xi_0)dk,
\end{equation}
which is further an integration
on a subset of $G_{n, r}$.
We claim the set  of this integration, namely
\begin{equation}
  \label{eq:int-set}
S:=\{\xi=
\exp(\frac{\pi}2 E)
k
\exp(\frac{\pi}2 E)
\cdot \xi_0; k\in K\},
\end{equation} 
 is given by
$$
S=\{\xi\in G_{n, r}; \xi 
\subset \eta_{1}\}
$$
where $\eta_1\in G_{n, n-r}$ is the element
$$
\eta_1:=\mathbb K^{r}\oplus 0\oplus  \mathbb K^{n-2r}
\subset
 \mathbb K^{r}\oplus \mathbb K^{r}
\oplus \mathbb K^{n-2r} =\mathbb K^n.
$$
Indeed let $\xi=
\exp(\frac{\pi}2 E)
k
\exp(\frac{\pi}2 E)\cdot \xi_0
$ be any element in the set
$S$ for some $k\in K$. Let $k=\diag(A, D)\in U(r)\times U(n-r)$
and write $D$  as a $2\times 2$-block matrix
under the decomposition of $\mathbb K^{n-r}
=\mathbb K^{r}\oplus \mathbb K^{n-2r}$,
$$D=\begin{bmatrix}
D_{11}& D_{12}\\
D_{21}& D_{22}
\end{bmatrix}.
$$
Using the formula   (\ref{eq:expX})
for the exponential we
find that 
$\xi$ is of the form 
$$
\xi=\{ D_{11} v \oplus  0 \oplus  D_{21}v; \,
v\in \mathbb K^{r}\} \subset \eta_1.
$$
Conversely, it is rather elementary
to see that any $r$-dimensional subspace
of $\eta_1$ is of the above form for some $D$. This
proves our claim. Therefore,
the integration
  (\ref{eq:int-form-sp}) is precisely the Radon transform of $\phi_{\m}$, namely
$$
\phi_{\m}(\xi_1)^2
=\mathcal R_{n-r, r}\phi_{\m}(\eta_1). 
$$
Now the function $\mathcal R_{n-r, r}\phi_{\m}$
is a $K$-invariant function on $G_{n, n-r}$,
and its dual Radon transform
evaluated at $\xi_0$,
$ \mathcal R_{n-r, r}^\ast (\mathcal R_{n-r, r}\phi_{\m}) (\xi_0)
= \mathcal R_{r, n-r} (\mathcal R_{n-r, r}\phi_{\m})(\xi_0)
$ is an integration of $\mathcal R_{n-r, r}\phi_{\m}$
over the set
$\{\eta; \eta\supset \xi_0\}$
which is a $K$-orbit of $\eta_1$. Namely the integrand
is constant and
$$\mathcal R_{r, n-r} (\mathcal R_{n-r, r}\phi_{\m})(\xi_0)
=(\mathcal R_{n-r, r}\phi_{\m})(\eta_1).
$$
Putting those together
we have
$$
\phi_{\m}(\xi_1)^2
=\mathcal R_{n-r, r}\phi_{\m}(\eta_1)
= \mathcal R_{n-r, r}^\ast
 \mathcal R_{n-r, r}\phi_{\m}(\xi_0).
$$
It follows from the main result in \cite{Grinberg-grs} (see Theorem 5.2 below)
that $ \mathcal R_{n-r, r}^t \mathcal R_{n-r, r}\phi_{\m}(\xi_0)$
is a non-zero constant of $\phi_{\m}(\xi_0)=1$
since $r\le n-r$. Thus $\phi_{\m}(\xi_1)\ne 0$.
 \end{proof}

 \begin{rema+} One may follow the argument in the
proof and find the constant
$\phi_{\m}(\xi_1)^2$ by using the result of Grinberg. It would
be interesting to evaluate
the constant for a general root system of positive multiplicities.
For the real projective space $G_{n, 1}=P^{n-1}(\mathbb R)$ the second claim (2)
is proved in \cite[ Chapter I, Lemma 4.9]{He2}.
 \end{rema+}

We prove now Theorem 4.3.
\begin{proof} We perform the spherical  transform
on the Bernstein-Sato formula for the cosine function in Theorem \ref{thm-1} with $\delta=2\nu +2$. Using
the self-adjoint property of the operator $\mathcal M_{2\nu}$ and
that (see \cite{Opdam-acta}),
$$
\mathcal M_{2\nu+2} \phi_{\m}
=\prod_{j=1}^r( (2\nu +2+\rho_1)^2 -(m_j+\rho_j)^2)\phi_{\m},
$$
we find that (suppressing the subindex $r$ in 
$c_{\nu, r}(\m)$)
$$
\frac {c_{\nu}(\m)}{N_\nu}
= \frac {N_{\nu+1}}{N_\nu}
\frac{\prod_{j=1}^r( (2\nu +2+\rho_1)^2 -(m_j+\rho_j)^2)}
{\prod_{j=1}^r(2\nu +2+a(j-1))(2\nu +1 +\iota +a(r-j))}
\frac 1{N_{\nu+1}}  c_{\nu+1}(\m).
$$
After a simplification we get
\begin{equation*}
\frac {c_{\nu}(\m)}{N_\nu}   =
\prod_{j=1}^r
(1-\frac
{\frac{m_j}2}
{\nu +1 +\frac a2(j-1)})
(1+\frac
{\frac{m_j}2}
{\nu +1 +b +\iota +a(r-1)  -\frac a2(j-1) }
)
\frac{c_{\nu+1}(\m)}{N_{\nu+1}}.
  \end{equation*}
Iterating the result produces furthermore
\begin{equation*}
\begin{split}
\frac {c_{\nu}(\m)}{N_\nu}  & =
\frac{c_{\nu+l+1}(\m)}{N_{\nu+l+1}}
\prod_{j=1}^r\prod_{k=0}^l
(1-\frac
{\frac{m_j}2}
{\nu +k+1  +\frac a2(j-1)})\\
&\quad 
(1+\frac
{\frac{m_j}2}
{\nu+k +1 +b +\iota +a(r-1)  -\frac a2(j-1) }
).
 \end{split}
  \end{equation*}
However $\frac{c_{\nu+l+1}(\m)}{N_{\nu+l+1}}\to \phi_{\m}(0)=1$, $l\to \infty$, according to Lemma 4.4.
Therefore,
\begin{equation*}
\begin{split}
\frac {c_{\nu}(\m)}{N_\nu}  & =
\prod_{j=1}^r\prod_{k=0}^\infty
(1-\frac
{\frac{m_j}2}
{\nu +k+1  +\frac a2(j-1)})\\
&\quad 
(1+\frac
{\frac{m_j}2}
{\nu+k +1 +b +\iota +a(r-1)  -\frac a2(j-1) }
),
 \end{split}
  \end{equation*}
which can also be written in terms of the Gamma function
(\cite[p.5]{Erdelyi-1})
$$
\prod_{j=1}^r 
\frac{ \Gamma(\nu +1 +\frac a2(j-1))
 \Gamma(\nu +1 +b +\iota+a(r-1)-\frac a2(j-1))
 }
{ \Gamma(\nu +1 +\frac a2(j-1) -\frac{m_j}2
)
\Gamma(\nu +1 +b +\iota+a(r-1)-\frac a2(j-1) +\frac{m_j}2).
}
$$
This proves our formula for the cosine function. The sine function
is done similarly.
\end{proof}

\section{Spectral symbols and range characterization
}

\subsection{Diagonalization of the transforms}
Theorem 4.3 applied to Lemma 3.5 gives then
\begin{prop+} 
Let $\mathbb K=\mathbb R, \bc$ or $ \mathbb H$,
and  $1\le r\le  n-1$.  
The eigenvalue
of $\mathcal C_{r, r}^{(\nu)}$ and  $\mathcal S_{r, r}^{(\nu)}$
are given by $c_{\nu, r}(\m)$ 
and $s_{\nu, r}(\m)$ in Theorem 4.3 with $\iota =a-1$, 
$2b=a(n-2r)$.
 \end{prop+}

To state our result on the spectral symbol of 
$\mathcal C_{r, r^\prime }^{(\nu)}$ and  $\mathcal S_{r, r^\prime}^{(\nu)}$
for different $r$ and $r^\prime$ we recall first
the following result of Grinberg \cite{Grinberg-grs}, reformulated
slightly differently here.

 \begin{theo+} Let $\mathbb K=\mathbb R, \bc$ or $ \mathbb H$,
and   $1\le r < r^\prime\le n-1$. 
Then the operator $\mathcal R_{r^\prime, r}
$ defines a
bounded operator from $L^2(G_{n, r})$ 
into $L^2(G_{n, r^\prime})$. 
The operator
$\mathcal R_{r^\prime, r}^\ast \mathcal R_{r^\prime, r}
$ is a diagonal
operator under the decomposition,
$$
\mathcal R_{r^\prime, r}^\ast
\mathcal R_{r^\prime, r} f
= c(\m)f, \quad f  \in V^{\m}, 
$$
with an explicit formula for the eigenvalue $c(\m)$.
The closure in $L^2(G_{n, r^\prime})$
of the image of the operator $\mathcal R_{r^\prime, r}$
on $L^2(G_{n, r})$ is
$$
\sum_{\m \in L_{r, r^\prime} }W^{\m} 
$$
where $L_{r, r^\prime}$ is the subset
of those $\m$ for which $ m_j =0$ if $j\ge \min\{r, r^\prime\}$.
\end{theo+}

Note that the $L^2$-bounded result was not stated in \cite{Grinberg-grs}. However
it follows directly from the explicit formula for the eigenvalue  $c(\m)$ there.
We remark also that the explicit formula for
$c(\m)$ found in \cite{Grinberg-grs} is 
under the condition $r\le r^\prime$ and $2r^\prime \le n$.
For general $r$ and $r^\prime$ the formula
for $c(\m)$ is given in \cite[Theorem 6.4]{Kakehi-jfa99}. 

Using Lemma 3.4, Proposition 5.1 and Theorem 5.2 we get
 \begin{coro+}
\label{coro3}
Suppose $\nu> 0$.  Let $1\le r,  r^\prime\le  n-1$. The eigenvalue of 
$\mathcal C_{r, r^\prime}^\ast 
\mathcal C_{r, r^\prime}$ and  $\mathcal S_{r, r^\prime}^\ast
\mathcal S_{r, r^\prime}
$
on the space $V^{\m}$ are given  respectively by
$$
c(\m)c_{\nu, r}(\m)c_{\nu, r^\prime}(\m), \qquad
s(\m)s_{\nu, r} (\m)s_{\nu, r^\prime}(\m),
$$
where 
$c_{\nu, r}$ is given in Proposition 5.1
and $c(\m)$  in Theorem 5.2
and  \cite{Kakehi-jfa99} .
In particular  $\mathcal C_{r, r^\prime}$ and 
$\mathcal S_{r, r^\prime}$ are bounded operators from
 $L^2(G_{n, r^\prime})$ to  $L^2(G_{n, r})$.
 \end{coro+}

\subsection{ Characterization of  the image 
of the transforms}
The following theorem follows immediately from Corollary 5.3 and
Lemma 4.5.

\begin{theo+}
 Let $1\le r,  r^\prime< n-1$ and $\nu >0$. 
 \begin{enumerate} 
\item 
Let  $\mathbb K=\br$. If $\nu\notin \frac{\mathbb Z}2$. Then
the  closures of images 
of $\mathcal C_{r^\prime, r}
$ and  $\mathcal S_{r^\prime, r}$, $L^2(G_{n, r}) \to
L^2(G_{n, r^\prime})$, 
are given by 
\begin{equation}\label{ima-1}
\sum_{\m\in L_{r, r^\prime}}W^{\m}.
\end{equation}
where $L_{r, r^\prime}$ is given in Theorem 5.2. In particular
the images are dense if $r^\prime \le r$.
 If $\nu\in\frac{\mathbb Z}2$
then the closures of their images in $L^2(G_{n, r^\prime})$ are given by
$$
\sum_{\m\in L_{r, r^\prime}\cap L_\nu }W^{\m},
$$
where $L_{\nu}$ is the subset of $\m$ such that
\begin{equation}
\label{L-nu-cond-1}
\frac{m_j}2 < \nu + 1+\frac 12 (j-1), \quad \text{if $\nu + 1+\frac 12 (j-1)$
is an integer}, \, \, j=1, \cdots, r.
\end{equation}
\item Let  $\mathbb K=\bc$ or $ \mathbb H$. If $\nu\notin  Z$. Then
the closure of the images are as in
(\ref{ima-1}).
If $\nu\in \mathbb Z$
then the 
closures are given by
$$
\sum_{\m\in L_{r, r^\prime}\cap L_\nu }W^{\m} ,
$$
where 
$L_{\nu}$ is the subset of $\m$ such that
\begin{equation}
\label{L-nu-cond-2}
\frac{m_j}2 < \nu + 1+\frac a2 (j-1), \quad  j=1, \cdots, r.
\end{equation}
 \end{enumerate} 
\end{theo+}

\section{The sine-transform 
as Knapp-Stein intertwining operator. Cosine transform and
branching of holomorphic representations}

In this section we give applications of our results
to the existence of the
Stein's complementary series and on the branching
rule of holomorphic representations.
\subsection{Stein's complementary series}

We fix $n=2r$ in this subsection. Let $GL_n=GL(n, \mathbb K)$
be the general linear group over $\mathbb K$ and $G=U(n, \mathbb K)$
as before. Consider the parabolic subgroup $P$ 
of $GL_n$ consisting
of block matrices of the form
$$
p=\begin{bmatrix} B&C\\
0&D
\end{bmatrix}
$$
where $B, C, D\in M_{r, r}(\mathbb K)$. 
 The Langlands decomposition of $P$ is $P=L N=MAN$,
with the nilpotent group $N$ 
and its opposite  $\bar N$  consisting   of upper
respectively lower triangular matrices 
\begin{equation}
  \label{eq:nil-iden}
n_X=\begin{bmatrix} I&X\\
0&I
\end{bmatrix}\in N, \quad n_Y=\begin{bmatrix} I&0\\
Y&I
\end{bmatrix}\in \bar N,
 \quad X, Y\in M_{r, r},
\end{equation}
both being identified with $M_{r, r}$,
the group  $L=MA$ consists of diagonal matrices
$$
p=\begin{bmatrix} B&0\\
0&D
\end{bmatrix}
$$
with  $B, D\in GL_r$ and $A$ being a one-dimensional
subgroup in the center of $L$.

 Let  $\delta_t$,  be    the one-dimensional representation 
$$
\delta_t(p)
= |\det_{\mathbb R}(BD^{-1})|^{t},
$$
of $P$, for $t\in \mathbb C$. 
The special case $\delta_{\frac r2}$
 corresponds to the 
determinant of adjoint representation of $P$ 
on $N$ (written generally as $e^{\rho}$
where $\rho$ is half-sum of positive roots
with respect to the one-dimensional  Lie
algebra of $A$).

Let $\Ind_{P}^{GL_n}(t)$ be the induced representation \cite[Chapter VII]{Kn-book} consisting 
of measurable functions $f$ on $G$ such  that
$$
f(gp)=(\delta_{t+\frac r2}(p)
)^{-1}f(g)
$$
and that $f|_G \in L^2(G)$.
When realized on $\bar N=M_{r, r}=
M_{r, r}(\mathbb K)$
the Knapp-Stein intertwining operator from
$\Ind_{P}^{GL_n}(t)$ to $Ind_{P}^{GL_n}(-t)$
is given by \cite{Stein-slnc}, 
$$
\mathcal I_t F(\bar n_X)=\int_{M_{r, r}} |\det_{\mathbb R}(X-Y)|^{(-r+ 2t )}
F(n_Y) dY.
$$

We will use however
the so-called compact realization
of the induced representation,  namely
the representation space
is the $L^2$-space $L^2(GL_n/{GL_n\cap L})$
on  the compact manifold $GL_n/{GL_n\cap L}$.
 Now $\mathcal X= GL_n/{GL_n\cap L}=G/K$ is precisely
the Grassmannian manifold $\mathcal X=G_{n, r}$.
The relation between the two
realizations is given by a change
of variables,
\begin{equation}
  \label{eq:k-to-n}
f(\xi)\in L^2(\mathcal X) \to  F(\bar n) =
\delta_{t+\frac r2}(p(\bar n))^{-1}
   f(\kappa(\bar n)\xi_0)  
\end{equation}
where $\bar n=k(\bar n)p(\bar n)$ is the Iwasawa $GL_n=GP$ 
decomposition of  $\bar n$ in  $GL_n$; see \cite[Chapter VII]{Kn-book}
(our $GL_n$ and $G$ correspond to $G$ respectively
$K$ there).

 \begin{lemm+} In the compact picture the Knapp-Stein 
(formal) intertwining
operator $\mathcal I_t$   is given by (up to a non-zero constant)
the sine transform $\mathcal S^{(\nu)}$,
$$
\mathcal J_t f(\xi)=
\mathcal S^{(\nu)} f(\xi)=
\int_{G_{n, r}}|\Sin (\xi, \eta)|^{2\nu} f(\eta) d\eta, 
$$
where
\begin{equation}
  \label{eq:nu-t}
\boxed{\nu=-\frac a2(r-2t)  }.
\end{equation}
It is well-defined on the space
$L^2(\mathcal X)$ for $\nu\ge 0$.
 \end{lemm+}
 \begin{proof} 
Let $S(\xi, \eta)$ denote temporarily
the kernel of the intertwining operator in the compact realization
on $L^2(\mathcal X
)$, namely let
$$
\mathcal J_t f(\xi)=\int_{\mathcal X}
S(\xi, \eta)f(\eta) d\eta.
$$
The kernel $S(\xi, \eta)$  is then uniquely
determined by $S(\xi_0, \eta) $ 
(recalling that
$\xi_0=\bbK^r \oplus  0\in 
\mathcal X$) by the transitivity
of $G=U(n, \mathbb K)$ on
$
\mathcal X
$ and by the intertwining property.
The evaluation 
$\mathcal J_t f(\xi_0)$ can also be computed 
via the intertwining operator $\mathcal I_t$ in the non-compact
picture above, 
\begin{equation}\label{J}
\int_{\mathcal X
} S(\xi_0, \eta )
f(\eta) d\eta=
\mathcal J_t f(\xi_0)
=\mathcal I_t F(0)
=\int_{M_{r, r}} |\det_{\mathbb R}(Y)|^{(-r+ 2t )}
F(n_Y) dY
\end{equation}
where the function $f$ on $\mathcal X
$ and 
$F$ on $\bar N$ is given by  (\ref{eq:k-to-n}). Performing
the  change
of variable $Y\in M_{r, r}\mapsto \eta:=k(n_Y)\xi_0\in 
\mathcal X
$
in the integral over $M_{r, r}$,  according to  (\ref{eq:k-to-n}),
 we find that
\begin{equation*}
 S(\xi_0, \eta )=|\det_{\mathbb R}(Y)|^{(-r+ 2t )}
\delta_{t+\frac r2}(p(n_Y))^{-1} 
(\text{Jac}_{Y\mapsto\eta})^{-1},
\end{equation*}
where $\text{Jac}_{Y\mapsto\eta}$ is the Jacobian of $Y\to \eta$, namely
$$
d\eta=\text{Jac}_{Y\mapsto\eta} dY.
$$

We will find the two quantities
$\delta_{t+\frac r2}(p(n_Y))$  and the Jacobian.
Consider therefore the $GL_n=GMAN$ decomposition
(note that $G$ is the maximal compact subgroup of $GL_n$)
of $n_Y\in \bar N= M_{r, r}(\mathbb K)$ (under
the identification (\ref{eq:nil-iden}))
\begin{equation}
\label{Iwas-n}
n_Y=\begin{bmatrix} I_r & 0\\
Y & I_r\\
\end{bmatrix}
=k( n_Y) p(n_Y),  \quad k(n_Y)\in G, \, p(n_Y)
=\begin{bmatrix} B&C\\
0&D\end{bmatrix}
\in P=MAN.
\end{equation}
This is a factorization of the
lower triangular matrix $n_Y$
as a product of a unitary
matrix $k(n_Y)$
and an upper triangular matrix $p(n_Y)$. 
Let the above element act on the base element $\xi_0\in \mathcal X$
via the defining action of $GL_n$ acting on subspaces of $\mathbb K^n$.
Since $p(n_Y)$ stabilize $\xi_0$
we see that the change of variables $Y\to \eta=k( n_Y)\xi_0\in \mathcal X$ is given by
\begin{equation}
\label{Y-eta}
\eta=k(n_Y)\xi_0=\begin{bmatrix} I_r &0 \\ Y &I
\end{bmatrix}
\xi_0=\{v \oplus Yv; v\in \mathbb K^r\}\in \mathcal X.
\end{equation}
An easy matrix computation then shows
also that
$$
\delta_{t+ \frac r2 }(p(n_Y))=\det_{\br}
(1+Y^\ast Y)^{t+\frac r2}.
$$
Indeed  taking 
the determinant of (\ref{Iwas-n}) we find that $
\det_{\mathbb R} B \det_{\mathbb R} D=1$; 
computing the upper left entry of the block matrix
 $n_Y^\ast n_Y$ we get
$B^\ast B=I +Y^\ast Y$ and 
$\delta_{s}(p(n_Y))=\det_{\mathbb R}(BD^{-1})^s
=\det_{\mathbb R}(B)^{2s}
=\det_{\mathbb R}(B^\ast B)^s=
\det_{\mathbb R}(I+Y^\ast Y)^s$ for any $s$.

The  measure $d\eta$, and equivalently the Jacobian, is given by
$$
d\eta=\det_{\mathbb R}
(I+Y^\ast Y)^{-r}dY, 
$$
by direct computations (or by using Proposition 3.3 (ii) in \cite{gz-radon}).
Now the sine $\Sin(\eta, \xi_0)$
of the angle between $\eta$ and $\xi_0$  is, by
(\ref{Y-eta}) and  Lemma 3.2,  given by
$$
|\Sin(\xi_0, \eta)|^{a}= \frac{|\det_{\mathbb R} (Y)|}
{\det_\br (1+Y^\ast Y)|^{\frac 12}}.
$$

It follows finally that 
$$
S(\xi_0, \eta )=|\Sin(\xi_0, \eta)|^{-a(r-2t)  }, 
$$
as claimed. This completes the proof.
 \end{proof}

\begin{rema+}Consider the simplest case when $GL_n=GL_2(\mathbb R)$.
The compact picture of the induced
representation in question is on the half unit circle $S^1/Z_2=\{
e^{i\theta}, 0\le \theta<\pi \}$ representing lines $\mathbb R(\cos\theta,
\sin\theta)$ in $\mathbb R^2$.
The intertwining operator in the non-compact picture
is an integral operator defined 
by the convolution by $|x|^{s}$ on $\mathbb R$.
The change of variables from non-compact to
the compact is then $x\to \mathbb R(1, x)$,
namely $x\to \mathbb R(\cos\theta, \sin\theta)$,
with $|\sin\theta| =\frac{|x|}{(1+x^2)^{\frac 12}}$.
The above lemma states that
the intertwining operator in the compact picture
is then the convolution
on the circle by $|\sin\theta|^s 
=(\frac 12|1-e^{2\theta}|)^s$, namely
it has the kernel $|e^{2\phi}-e^{2\theta}|^s$,
which is well-known.
  \end{rema+}
By using Lemma 3.5 and Theorem 4.3 we obtain then
the eigenvalues of the Knapp-Stein intertwining
operator $\mathcal J_t$ on $L^2(G_{n, r})$ under the decomposition
(2.4), and we can determine the range of $t$ for which
the eigenvalues are all non-negative (or all non-positive).
This  gives   an independent proof of
the existence proof Stein's complementary series
and with explicit formula for  the corresponding  unitary structure.

\begin{theo+}\label{thm4}
 Let $\nu=-\frac a2(r-2t)
$.  The Knapp-Stein intertwining
operator $\mathcal J_t$, 
 is well-defined for $\nu \ge 0$
on
the algebraic span of the subspace $V^{\m}$. It intertwines
the actions of $\mathfrak {gl}(n, \mathbb K)$
by the induced representations 
$\Ind_{P}^{GL_n}(t)$ and $Ind_{P}^{GL_n}(-t)$, and has
meromorphic continuation in $\nu$.
 Its eigenvalues are given  by
\begin{equation}\label{stein-com-ser}
N_{\nu}^\prime\prod_{j=1}^r\frac{(\frac a2(r+1-j -2t))_{\frac {m_j}2} }
{(\frac a2(r+j-1+2t))_{\frac {m_j}2}}.
\end{equation}
The sesqui-linear
form
$$
(f, g)_\nu:=
\frac 1{N_\nu^\prime} (\mathcal J_t f, g)
$$
for
$$
0<t<\frac 12,
$$
is well-defined, and  is a $\mathfrak {gl}(n, \mathbb K)$-invariant
positive definite Hermitian
inner product.
 The  completion of the pre-Hilbert space is a unitary
 representation of  $GL(n, \mathbb K)$.
 \end{theo+}

For  more details on unitarity and composition series
of the whole family of the induced representations
see also \cite{Howe-Lee-jfa},
 \cite{Sahi-crelle} and \cite{gkz-ma}.

\subsection{Branching
of holomorphic representations.}

Finally we give
an application of our result to
the  branching of holomorphic representations
on compact Hermitian symmetric spaces. The non-compact
case has been studied intensively;  see \cite{gkz-bere-rbsd} 
and references therein. The result below
can be deduced from the general theory 
\cite{kobayashi-06}
of (discrete) branching
of highest weight representations.
 We only indicate here
a concrete approach using the cosine transform.
To keep the presentation of paper rather explicit 
we will only treat the case  when the
Hermitian symmetric spaces is the complex Grassmannian
manifold with real form being the real or quaternionic
Grassmannians. 

We fix 
the complex Grassmannian manifold
 $\mathcal X_1$ of $r_1$-dimensional
complex subspace in $\mathbb C^{2r_1 +b_1}$,
$r_1\ge 1, b_1\ge 0$,
namely
$$
\mathcal X_1=U(2r_1+b_1)
/U(r_1)\times U(r_1+b_1).$$
$\mathcal X_1$ is equipped with
 the $U(2r_1+b_1)$-invariant Hermitian metric.
Consider the real Grassmannian $\mathcal X$,
$$\mathcal X=G/K=O(2r+b)/O(r)\times O(r+b),
\quad \, r=r_1
$$
or
the quaternionic Grassmannian
 $$\mathcal X=G/K=Sp(r +b)/
{Sp(r) \times Sp(b)}, \qquad r=\frac{r_1}2, \qquad
b=\frac{b_1}2
$$
 when $r_1$ and $b_1$ are even integers.
Then $\mathcal X$ can be realized as a
 totally geodesic real form of the Hermitian
symmetric space $\mathcal X_1$.
We realize
the space $V_1=M_{ r_1+b_1, r_1}(\mathbb C)$
as a dense subset  of $\mathcal X_1$ by the identification
$$z\in V_1 \mapsto
\{y \oplus  zy; y\in \bc^{r_1}\}\in \mathcal X_1;$$
 similarly
$V=M_{ r+b, r}(\mathbb K)$ 
can be realized as a dense subset of $\mathcal X$
via
$$z\in V \mapsto
\{y \oplus  zy; y\in \bbK^{r}\}\in \mathcal X.
$$
For any positive integer $\alpha $ there is
a corresponding weighted Bergman space
on the compact space $\mathcal X_1$, denoted
by $ H_{\alpha}(\mathcal X_1)$,  with the reproducing
kernel $\det (I+w^\ast z)^{\alpha}$, 
 $z, w\in V\subset \mathcal X_1$, which
forms an irreducible unitary representation
of $U(2r_1+b_1)$. The elements
in the space $ H_{\alpha}(\mathcal X_1)$ will be
identified with holomorphic polynomials on $V_1$.
Consider $ H_{\alpha}(\mathcal X_1)$
as a unitary representation of $G\subset U(2r_1+b_1)
$ and we
will find the explicit irreducible decomposition.

We can first realize
the  $L^2(\mathcal X)$ as a space of functions on the subset $V$. Indeed the 
$L^2(\mathcal X)$ is $G$-equivalent to
$$
L^2(V, \det(I+x^\ast x)^{-\frac{2r_1+b_1}2} dm(x))
$$
where $dm(x)$ is the Lebesgue measure.

The restriction mapping
$$
\mathcal T:  H_{\alpha}(\mathcal X_1) \to 
C^\infty(\mathcal X),  f\mapsto f(x) \det(I+x^\ast x)^{-\frac \alpha 2}
$$
defines then an intertwining map, where $G$ acts on $C^\infty(\mathcal X)$
by the defining action. The operator
$\mathcal T\mathcal T^\ast$ 
on $L^2(\mathcal X)=L^2(V, \det(I+x^\ast x)^{-\frac{2r_1+b_1}2} dm(x))$ 
will be called the Berezin transform and
is of the form
$$
\mathcal T
\mathcal T^\ast f(x)= \int_{V}
\det(I+x^\ast x)^{-\frac \alpha 2}
\det(I+x^\ast y)^{ \alpha } f(y)
\det(I+y^\ast y)^{-\frac \alpha 2}
\det(I+y^\ast y)^{-\frac{2r_1+b_1}2} dm(y);
$$
see \cite{gkz-bere-rbsd}.
This is just the cosine transform, and
its eigenvalues on the space $L^2(\mathcal X)$
implies  the branching rule of the
representation of $U(2r_1+b_1)$
under $G$. The first part of the next Proposition
follows by similar computation as in the previous subsection,
 the second part follows from Theorem 5.4 and some abstract
arguments.

\begin{prop+}\label{branching}
\begin{enumerate}
\item The Berezin transform $\mathcal T\mathcal T^\ast $
on $\mathcal X=G/K=G_{n, r}(\mathbb K)$ is the cosine transform $\mathcal 
C_{r, r}^{(\nu)}$
with $\nu =\frac{\alpha}2$ for $\bbK =\br $ and $\nu =\alpha$ for $\bbK =\mathbb H $.
\item  The representation  $ \mathcal H_{\alpha}(\mathcal X_1)$
is decomposed under $G$ with multiplicity free 
as 
 $$
 \mathcal H_{\alpha}(\mathcal X_1)=\sum_{\m \in L_{\nu}}V^{\m},
$$
where $L_\nu$ is given in
(\ref{L-nu-cond-1})  and (\ref{L-nu-cond-2}).
\end{enumerate}
\end{prop+}

\def\cprime{$'$} \newcommand{\noopsort}[1]{} \newcommand{\printfirst}[2]{#1}
  \newcommand{\singleletter}[1]{#1} \newcommand{\switchargs}[2]{#2#1}
  \def\cprime{$'$} \def\cprime{$'$}
\providecommand{\bysame}{\leavevmode\hbox to3em{\hrulefill}\thinspace}
\providecommand{\MR}{\relax\ifhmode\unskip\space\fi MR }
\providecommand{\MRhref}[2]{%
  \href{http://www.ams.org/mathscinet-getitem?mr=#1}{#2}
}
\providecommand{\href}[2]{#2}

\end{document}